\hfill {\bf IC/97/110} 
\baselineskip=16pt

\font\twelvebf=cmbx12

\vskip 15mm
\def\[{[\![}
\def\]{]\!]}
\def\l{\langle}
\def\r{\rangle}

\def\C{{\bf C}}
\def\Ch{{\C[[h]]}}

\def\Q{{\bf Q}}
\def\s{\smallskip}
\def\e{\eqno}
\def\n{\noindent}
\def\Uq{U_q[osp(2n+1/2m)]}
\def\U{U[osp(2n+1/2m)]}
\def\q{{\bar q}}
\def\k{{\bar k}}
\def\L{{\bar L}}

\def\a{\alpha}

\hfill    {\bf ICTP Preprint IC/97/110}

\vskip 30mm

\n
{\twelvebf A q-deformation of the parastatistics and
an alternative to the Chevalley  \hfill\break description of
U$_{\bf q}$[osp(2n+1/2m)] }

\vskip 32pt
\noindent
T.D. Palev\footnote*{Permanent address: Institute for Nuclear
Research and Nuclear Energy, 1784 Sofia, Bulgaria; E-mail:
tpalev@inrne.acad.bg}

\noindent
International Centre for Theoretical Physics, 34100 Trieste, Italy
\vskip 12pt

\vskip 4cm

\vskip 48pt
\noindent
{\bf Abstract.} The paper contains essentially two new results.
Physically, a deformation of the parastatistics in a sense of
quantum groups is carried out. Mathematically, an alternative to
the Chevalley description of the quantum orthosymplectic
superalgebra $\Uq$ in terms of $m$ pairs of deformed parabosons
and $n$ pairs of deformed parafermions is outlined.

\vfill\eject
\leftskip 0pt
\vskip 12pt

\n
{\bf 1. Introduction}

\bigskip
In this paper we give an alternative to the Chevalley
definition of the quantum superalgebra $\Uq$ in terms of
generators and relations (see Eqs. (30)). We generalize to the quantum
case a result we have recently obtained,$^1$  namely that the
universal enveloping algebra $\U$ of the orthosymplectic Lie
superalsebra $osp(2n+1/2m)$ is an associative unital algebra with
generators, called Green generators (operators),
$$
a_1^\pm,a_2^\pm,\ldots,a_{m-1}^\pm,a_{m}^\pm,a_{m+1}^\pm,
\ldots,a_{m+n}^\pm \equiv a_N^\pm,\e(1)
$$
and relations
$$
[[a_{N-1}^\eta,a_{N}^{\eta}],a_{N}^\eta]=0,\quad
\[\[a_i^\eta,a_j^{-\eta}\],a_k^\eta\]=
2\eta^{\l k \r} \delta_{jk}a_i^\eta, \quad \forall \; |i-j|\leq 1,
\quad \eta =\pm , \e(2) 
$$
where
$$
{\rm deg}(a_i^\pm)\equiv \langle i \rangle = \cases 
{\bar 1, & for $i\le m$ \cr
\bar 0, &  for $i> m$. \cr} \e(3)
$$
Here and throughout $\[a,b\]=ab-(-1)^{deg(a)deg(b)}ba$,
$[a,b]=ab-ba$, $\{a,b\}=ab+ba.$

The motivation for the present work stems from the observation
that the Green operators provide a description
of $osp(2n+1/2m)$ via generators, which, contrary to the
Chevalley elements, have a direct physical significance. 
As it was shown in Ref. 1
$a_1^\pm,a_2^\pm,\ldots,a_{m-1}^\pm,a_{m}^\pm$ 
(resp. $a_{m+1}^\pm,\ldots,a_{m+n}^\pm \equiv a_N^\pm)$
are para-Bose (pB) (resp. para-Fermi (pF)) operators. 
These operators were introduced in the quantum field theory as a possible
generalization of the statistics of the tensor (resp. spinor)
fields.$^2$ Therefore what we are actually doing here is a simultaneous
deformation of the pB and the pF operators (2) in the sense
of quantum groups.$^{3, 4}$

The fact that $n$ pairs of pF creation and annihilation operators
(CAOs) generate the Lie algebra so(2n+1) was first observed in
Refs. 5 and 6. It took some time to incorporate the para-Bose statistics
into an algebraic structure: $a_{1}^\pm,\ldots,a_{m}^\pm $ are
odd elements, generating a Lie superalgebra $^7$ isomorphic to
$osp(1/2m)$.$^8$ It is usually assumed that the pB operators
commute with pF operators. Other possibilities were also 
investigated.$^9$ In Ref. 10 it was indicated that the relations
between the pB and pF operators can be selected in such a way that
they generate $osp(2n+1/2m)$.

The identification of the parastatistics with a well known
algebraic structure has far going consequences. Firstly, it
indicates that the representation theory of $n$ pairs of pF
operators (or of $m$ pairs of pB operators) is completely
equivalent to the representation theory of $so(2n+1)$ (resp. of
$osp(1/2m)$). On this way one may enlarge considerably $^{11,
12}$ the class of the known representations, those corresponding
to a fixed order of the statistics. In particular, since the
(complex) Lie algebra $so(2n+1)$ has infinite-dimensional
representations, so do the pF operators.  Similarly $osp(1/2m)$
has finite-dimensional representations (for instance the defining
one) and therefore the pB operators have also such
representations. Secondly, it provides a natural background for
further generalizations of the quantum statistics.  In order to
give a hint of where this possibility comes from consider in the
frame of the quantum field theory a  field $\Psi(x)$. In the
momentum space the translation invariance of the field is
expressed as a commutator between the energy-momentum $P^m, \;
m=0,1,2,3$ and the CAOs $a_i^{\pm}$ of the field: 
$$ 
[P^m, a_i^\pm]=\pm k_i^m a_i^\pm ,\e(4) 
$$ 
where the index $i$ replaces
all (continuous and discrete) indices of the field and 
$$
P^m=\sum_i k_i^m H_i.   \e(5) 
$$ 
In the case of the pF
statistics $H_i={1\over 2}[a_i^+,a_i^-]$, whereas for pB fields
$H_i={1\over 2}\{a_i^+,a_i^-\}$.  In a unified form $H_i={1\over
2}\[a_i^+,a_i^-\]$ with the pF considered as even elements and
the pB - as odd. To quantize the field means, loosely speaking, to
find solutions of Eqs. (4) and (5), where the unknowns are the
CAOs $a_i^\pm$.

{\it The first opportunity} for further generalizations is based
on the observation that both $so(2n+1)$ and $osp(1/2m)$ belong to
the class {\it B} superalgebras in the classification of
Kac.$^{13}$ Therefore it is natural to try to satisfy the
quantization equations (4) and (5) with CAOs, generating
superalgebras from the classes {\it A}, {\it C} or {\it D}
$^{13}$ or generating other superalgebras from the class {\it B}.
It turns out this is possible indeed. Examples of this kind, notably
the $A-$statistics, related to the completion and the central
extension of $sl_\infty$ were studied in Refs. 14 and 15 (see
also {\it Example 2} in Ref. 16 and the the other references in
this paper). The Wigner quantum systems (WQSs), introduced in
Refs. 17  and 18, are also examples of this kind, however in the
frame of a noncanonical quantum mechanics.  Some of these systems
possess quite unconventional physical features, properties which
cannot be achieved in the frame of the quantum mechanics. The
$(n+1)-$particle WQS, based on $sl(1/3n) \in A$,$^{19}$ exhibits
a quark like structure: the composite system occupies a small
volume around the centre of mass and within it the geometry is
noncommutative.  The underlying statistics is a Haldane exclusion
statistics,$^{20}$ a subject of considerable interest in
condensed matter physics.  The $osp(3/2) \in B$ WQS, studied in
Ref. 21, leads to a picture where two spinless point particles,
curling around each other, produce an orbital (internal angular)
momentum $1/2$.

{\it The second opportunity} for generalization of the statistics
is based on deformations of the relations (2).  Assume only for
simplicity that in (5) $i=1,\ldots,n$  (the consideration remain
valid for $n=\infty$).  Then both in the pF and the pB case $H_i$
are elements from the Cartan subalgebra $H$ of $so(2n+1)$ and
$osp(1/2n)$, respectively.  The CAOs are  root vectors of these
(super)algebras (see Ref. 15 for more detailed discussions).  The
important point now comes from the observation that the
commutation relations between the Cartan elements and the root
vectors, in particular the quantization equations (4) and (5),
remain unaltered upon $q$-deformations.  Therefore one can
satisfy the quantization equations (4) and (5) also with deformed
pF (resp. pB) operators.  Certainly in this case the relations
$H_i={1\over 2}\[a_i^+,a_i^-\]$ cannot be preserved anymore. One
has to postulate the expression (5), introducing the additional
(Cartan) generators $H_i$ similarly as in the case of a deformed
harmonic oscillator, where one is forced to introduce also number
operators.$^{22,23,24}$

The conclusion is that the deformed pF operators $\{a_i^\pm,H_i
|i=1,\ldots,n\}$, being solutions of the quantization equations
(4) and (5), enlarge the class of the possible statistics.
It turns out these are the  operators, which provide an
alternative to the Chevalley description of $U_q[so(2n+1)]$. This
was shown in Ref. 25. A similar problem for $U_q[osp(1/2m)]$,
corresponding to a $q-$deformation of the pB operators, was first
carried out for $m=1^{26}$ and then for any $m$.$^{27,28,29}$
Here we generalize the results for any $U_q[osp(2n+1/2m)]$,
$n,m>1$, namely when both para-Bose and para-Fermi operators are
involved.  This amounts to a simultaneous deformation of the
parabosons and the parafermions as one single supermultiplet.

In Sect. 2, after recalling the definition of $\Uq$, we
introduce the deformed Green generators (13) and derive the
relations (30) they satisfy. In Sect 3 we solve the invere
problem.  We express the Chevalley elements via the Green
generators and show that the relations among the Chevalley
generators follow from the properties of the Green operators.
This leads to the conclusion that the deformed Green operators
provide an alternative description of $\Uq$.

\bigskip
Throughout the paper we use the notation
(some of them standard):

$\C$ - all complex numbers;

$\C[[h]]$ - the ring of all formal power series in $h$ over $\C$;

$q=e^h \in \Ch,  \quad {\bar q}=q^{-1} $;

$\[a,b\]=ab-(-1)^{deg(a)deg(b)}ba, \quad
[a,b]=ab-ba, \quad \{a,b\}=ab+ba;$

$\[a,b\]_x=ab-(-1)^{deg(a)deg(b)}xba, \quad
[a,b]_x=ab-xba, \quad \{a,b\}_x=ab+xba;$

$$
{\rm deg}(a_i^\pm)\equiv \langle i \rangle = \cases 
{\bar 1, & for $i\le m$ \cr
\bar 0, &  for $i> m$. \cr} \hfill
$$

$$
q_i=q^{(-1)^{\l i+1 \r}},\;\; i.e.\;\;
q_i={\bar q},\;i<m, \quad q_i=q,\;i\ge m.
$$

\bigskip
For a convenience of further references we list here some
deformed identities, which will be often used
($Id(3)$ follows from $Id(2)$).

\smallskip
\noindent
$Id(1)$: If $[a,c]=0$,
then $(x+x^{-1})[b,[a,[b,c]_x]_x]=[a,[b,[b,c]_x]_{x^{-1}}]_{x^2}-
[[b,[b,a]_x]_{x^{-1}},c]_{x^2}$;

\smallskip
\noindent
$Id(2)$: If $B$ or $C$ is an even element, then for any 
values of $x,y,z,t,r,s$ subject to the relations\hfill\break
${}  \quad \quad \quad x=zs,\;\; y=zr,\;\; t=zsr,$
$$
\[A,[B,C]_x\]_y=\[\[A,B\]_z,C\]_t +
(-1)^{deg(A)deg(B)}z\[B,\[A,C\]_r\]_s;
$$

\smallskip
\noindent
$Id(3)$: If  $C$ is an even element and $[A,C]=0$, then
$$
\[A,[B,C]_x\]_y=[\[A,B\]_y,C]_x.
$$

\vskip 24pt

\n
{\bf 2. Deformed  Green Generators and their relations}

\bigskip
The $q$-deformed superalgebra $U_q[osp(2n+1/2m)]$, a Hopf superalgebra,
is by now a classical concept. See, for instance, Refs. 30-33, 
where all Hopf algebra operations are explicitly given. Here,
following Ref. 33, we write only the algebra operations.

Let
$(\a_{ij}), \; i,j=1,\ldots,m+n=N$ be an $N\times N$ symmetric Cartan
matrix chosen as:
$$
(\a_{ij})= (-1)^{\l j\r}
\delta_{i+1,j}+(-1)^{\l i \r}
\delta_{i,j+1} -[(-1)^{\l j+1
\r} + (-1)^{\l j \r}]\delta_{ij}+
\delta_{i,m+n}\delta_{j,m+n}. \eqno(6)
$$
For instance the Cartan matrix, corresponding to $m=n=4$
is $ 8 \times 8$ dimensional matrix:

\bigskip
$$
(\a_{ij})=\pmatrix{
    2&-1& 0& 0& 0& 0& 0& 0\cr
   -1& 2&-1& 0& 0& 0& 0& 0\cr
    0&-1& 2&-1& 0& 0& 0& 0\cr
    0& 0&-1& 0& 1& 0& 0& 0\cr
    0& 0& 0& 1&-2& 1& 0& 0\cr
    0& 0& 0& 0& 1&-2& 1& 0\cr
    0& 0& 0& 0& 0& 1&-2& 1\cr
    0& 0& 0& 0& 0& 0& 1&-1\cr
}. \eqno(7)
$$

\bigskip\n
{\bf Definition.}{\it
$U_q[osp(2n+1/2m)]$ is a Hopf algebra, which is a
topologically free module over $\C[[h]]$ (complete in $h$-adic
topology), with (Chevalley) generators $h_i,\;e_i,\;f_i,$
$i=1,\ldots,N$ and 
\smallskip
\noindent
1. Cartan-Kac relations:
$$
\eqalignno{
& [h_i,h_j]=0, & (8a) \cr
& [h_i,e_j]=a_{ij}e_j,\quad
  [h_i,f_j]=-a_{ij}f_j, & (8b) \cr
& \[e_i,f_j\]=\delta_{ij}{{k_i-{\bar k}_i}\over{q-{\bar q}}}, 
  \quad  k_i=q^{h_i},\;
   {\bar k}_i = k_i^{-1}=q^{-h_i},& (8c) \cr
}
$$
2. $e-$Serre relations
$$
\eqalignno{
& \[e_i,e_j\]=0, \quad \vert i-j \vert \neq 1, & (9a) \cr
& [e_i,[e_i,e_{i \pm 1}]_{\bar q}]_q \equiv
  [e_i,[e_i,e_{i \pm 1}]_q]_{\bar q}=0, \quad
  i\neq m,\;\; i\neq N & (9b) \cr 
& \{[e_{m},e_{m-1}]_{q},[e_m,e_{m+1}]_{\bar q}\}
\equiv \{[e_{m},e_{m-1}]_{\q},[e_m,e_{m+1}]_{q}\}=0, & (9c)\cr
& [e_N,[e_N,[e_N,e_{N-1}]_{\bar q}]]_q \equiv
[e_N,[e_N,[e_N,e_{N-1}]_q]]_{\bar q}=0, & (9d)\cr 
}
$$
3. $f-$Serre relations
$$
\eqalignno{
& \[f_i,f_j\]=0, \quad \vert i-j \vert \neq 1, & (10a) \cr
& [f_i,[f_i,f_{i \pm 1}]_{\bar q}]_q \equiv
  [f_i,[f_i,f_{i \pm 1}]_q]_{\bar q}=0, \quad
  i\neq m,\;\; i\neq N & (10b)\cr 
& \{[f_{m},f_{m-1}]_{q},[f_m,f_{m+1}]_{\bar q}\}
  \equiv \{[f_{m},f_{m-1}]_{\q},[f_m,f_{m+1}]_{q}\}=0, & (10c) \cr
& [f_N,[f_N,[f_N,f_{N-1}]_{\bar q}]]_q \equiv
  [f_N,[f_N,[f_N,f_{N-1}]_q]]_{\bar q}=0. & (10d) \cr 
}
$$
The grading on $U_q[osp(2n+1/2m)]$ is induced from:
$$
deg(h_j)=\bar 0,\; \forall j,\quad
deg(e_m)=deg(f_m)=\bar 1,
\; deg(e_i)= deg(f_i)=\bar 0  \;\;
for \;\; i\neq m. \e(11)
$$
}

The (9c) and (10c) relations are the additional Serre
relations,$^{33, 34, 35}$ which were initially omitted.
We do not write the other Hopf algebra maps ($\Delta,\;
\varepsilon,\; S$) $^{33}$ since we will not use them. They are
certainly also a part of the definition of $\Uq$. 

From (8a,b) and the definition of $k_i$ one derives:
$$
\eqalignno{
&  k_ik_i^{-1}=k_i^{-1}k_i=1, \quad  k_ik_j=k_jk_i,& (12a)\cr
&  k_ie_j=q^{\alpha_{ij}}e_jk_i,\quad k_if_j=q^{-\alpha_{ij}}f_jk_i.&(12b)\cr
}
$$

Introduce the following $3N$ elements in  $U_q[osp(2n+1/2m)]$
($i=1,\ldots,N-1$):
$$
\eqalignno{
&  a_i^-=(-1)^{(m-i)\l i \r}\sqrt{2}[e_i,[e_{i+1},
  [\ldots,[e_{N-2},[e_{N-1},e_N]_{q_{N-1}}]_{q_{N-2}}
  \ldots ]_{q_{i+2}}]_{q_{i+1}}]_{q_{i}},\quad 
   a_N^-=\sqrt{2}e_N, & (13a)   \cr
&  a_i^+=(-1)^{N-i+1}\sqrt{2}
   [[[\ldots [f_N,f_{N-1}]_{{\bar q}_{N-1}},f_{N-2}]_{{\bar q}_{N-2}} 
   \ldots]_{{\bar q}_{i+2}},f_{i+1}]_{{\bar q}_{i+1}},f_{i}]_{{\bar q}_{i}}
   , \quad a_N^+=-\sqrt{2}f_N,  &(13b) \cr  
&  H_i=h_i+h_{i+1}+\ldots +h_N \;({\rm including} \;i=N),& (13c)  \cr
}
$$
We refer to the above operators as to (deformed) Green generators
since in the nondeformed case they coincide with the Green generators
of $U_q[osp(2n+1/2m)]$ (see Eq. (36) in Ref. 1). Therefore 
$a_1^\pm,a_2^\pm,\ldots,a_{m-1}^\pm,a_{m}^\pm$ 
(resp. $a_{m+1}^\pm,\ldots,a_{m+n}^\pm \equiv a_N^\pm)$
can be viewed as deformed para-Bose (pB) (resp. deformed
para-Fermi (pF)) operators. Our aim is to show that $\Uq$
can be described entirely via the generators (13).

One can write (13a) and (13b) as 
$$
\eqalignno{
& a_i^-=(-1)^{(m-i)\l i \r +(m-j)\l j \r}[e_i,[e_{i+1},[\ldots
[e_{j-2}, [e_{j-1},a_j^-]_{q_{j-1}}]_{q_{j-2}}\ldots ]_{q_{i+1}}]_{q_i},\quad
		 i<j< N,& (14a) \cr
& a_i^+=(-1)^{i+j}
     [[\ldots[[a_j^+,f_{j-1}]_{{\bar q}_{j-1}},f_{j-2}]_{{\bar q}_{j-2}}
     \ldots]_{{\bar q}_{i+2}},f_{i+1}]_{{\bar q}_{i+1}},
	 f_i]_{{\bar q}_{i}},\quad i<j< N.  & (14b)\cr
}
$$		  
Taking into account (9a), (10a) and applying repeatedly  $Id(3)$
one rewrites the Green generators also as
$(i<j< N)$ 
$$
\eqalignno{
& a_i^-=(-1)^{(m-i)\l i \r +(m-j)\l j \r}[[e_i,[e_{i+1},[e_{i+2},\ldots
        [e_{j-3},[e_{j-2}, e_{j-1}]_{q_{j-2}}]_{q_{j-3}}\ldots 
       ]_{q_{i+2}} ]_{q_{i+1}}]_{q_i},a_j^-]_{q_{j-1}}, & (15a)\cr
& a_i^+=(-1)^{i+j}[a_j^+,
     [\ldots [[f_{j-1},f_{j-2}]_{{\bar q}_{j-2}},f_{j-3}]_{\q_{j-3}}
     \ldots]_{{\bar q}_{i+2}},f_{i+1}]_{{\bar q}_{i+1}},
	 f_i]_{{\bar q}_{i}}]_{\q_{j-1}}  .& (15b)\cr 
}
$$

The next proposition plays an important role in several intermediate
calculations.

\smallskip\n
{\bf Proposition 1.} {\it The following "mixed" relations between the
Chevalley and the Green generators take place:}
$$
\eqalignno{
 & \[e_i,a_j^+\]=-\delta_{ij}(-1)^{\l i+1 \r}k_ia_{i+1}^+,
       \quad  i\ne N, & (16)\cr
 &&\cr
 & \[a_j^-,f_i\]=\delta_{ij}a_{i+1}^-{\bar k}_i,\quad i\ne N, & (17)\cr
 &&\cr
 & \[e_i,a_j^-\]=0, \;\;if\; i<j-1 \; or \; i>j,\quad i\ne N,& (18a)\cr
 & \[e_i,a_{i+1}^-\]_{q_i}=(-1)^{\l i+1 \r}a_i^-, 
       \quad i\ne N,   &  (18b) \cr
 & \[e_i,a_{i}^-\]_{{\bar q}_{i-1}}=0,\quad i\ne N, & (18c) \cr
 &&\cr
 & \[a_j^+,f_i\]=0, \;\;if\; i<j-1 \; or \; i>j,\quad i\ne N,& (19a)\cr
 & \[a_{i+1}^+,f_i\]_{{\bar q}_i}=-a_i^+,\quad i\ne N. & (19b) \cr
 &  \[a_i^+,f_i\]_{q_{i-1}}=0, \quad i\ne N. & (19c) \cr
}
$$

\n
{\it Proof.} We stress on some of the intermediate steps in the
proof. 

\s\n
1. Begin with (16).

\n
(i) Let $i<j$. Then from (13b) and (8c) one immediately has 
$\[e_i,a_j^+\]=0$.

\smallskip\n
(ii) Let $i>j$. From (14b)
$\[e_i,a_j^+\]\sim \[e_i,[[...[[[a_{i+1}^+,f_i]_{\q_i},f_{i-1}]_{\q_{i-1}},
f_{i-2}]_{\q_{i-2}}\ldots]_{\q_{j-1}},f_j]_{\q_{j}} \]  $
(applying repeatedly $Id(3)$ and (8c))
$=[\ldots[A,f_{i-2}]_{\q_{i-2}}\ldots]_{\q_{j-1}},f_j]_{\q_{j}}$,
where $A=[\[e_i,[a_{i+1}^+,f_i]_{\q_{i}}\],f_{i-1}]_{\q_{i-1}}$
(from (i) and $Id(3)$)
=$[[a_{i+1}^+,\[e_i,f_i\] ]_{\q_i},f_{i-1}]_{\q_{i-1}}\sim
[a_{i+1}^+k_i,f_{i-1}]_{\q_{i-1}}=
\q_{i-1}[a_{i+1}^+,f_{i-1}]k_i=0 $, since evidently
$f_{i-1}$ commutes with $a_{i+1}^+$ (see (13b). Hence 
$$
\[e_i,a_j^+\]=0 \quad for\;\; i>j. \e(20) 
$$

\smallskip\n
(iii) Let $i=j$. 
$
\[e_i,a_i^+\]=-\[e_i,[a_{i+1}^+,f_i]_{\q_i}\]$
(from (i) and $Id(3)$)$=-[a_{i+1}^+,\[e_i,f_i\]]_{\q_i}=
-[a_{i+1}^+,{{k_i-\k_i}\over {q-\q}}]_{\q_i}=
-(-1)^{\l i+1 \r}\q_i a_{i+1}^+k_i=
-(-1)^{\l i+1 \r}k_i a_{i+1}^+.$ The unification of
(i)-(iii) yields (16).

\s\n 
2. Eq. (17) is proved in a similar way. 

\s\n
3. We pass to prove (18a).

\smallskip\n
(i) The case $i<j-1$ is evident.

\smallskip\n
(ii) Take $i=m>j$.  Note first that according to (9) and $Id(3)$
$
\[e_m,[e_{m-1},[e_m,e_{m+1}]_q]_{\q}\]$ \hfill\break
$= \[e_m,[[e_{m-1},e_m]_{\q},e_{m+1}]_q\]$ (using $Id(2)$) 
$=\{ [e_m,e_{m-1}]_q,[e_m,e_{m+1}]_{\q}\}-
[e_{m+1},\{e_m,[e_{m-1},e_m]_{\q}\}_q]  =0$,
\hfill\break
according to (9c) and (9a), i.e.,
$$
B\equiv\[e_m,[e_{m-1},[e_m,e_{m+1}]_q]_{\q}\]=0. \e(21)
$$
If $m=N-1$, then $\[e_m,a_{m-1}^-\]=B=0$. Let $m<N-1.$ From
(15a)  the and $Id(3)$ $\[e_m,a_{m-1}^-\]\sim [B,a_{m+2}^-]_q=0.$

\n
(iii) Let $i\ne m >j$. From (15b) $a_{i-1}^- \sim 
[[e_{i-1},[e_i,e_{i+1}]_{q_i}]_{q_{i-1}},a_{i+2}^-]_{q_{i+1}}$.
Then $Id(3)$ yields
$[e_i,a_{i-1}^-]\sim [z,a_{i+2}]_{q_{i+1}}$ with
$z=[e_i,[e_{i-1},[e_i,e_{i+1}]_{q_i}]_{q_i}]$. 
If $i>m$, then using $Id(1)$, 
$z=[e_i,[e_{i-1},[e_i,e_{i+1}]_{q}]_{q}] \hfill\break
\sim [e_{i-1},[e_{i},[e_i,e_{i+1}]_{q}]_{\q}]_{q^2}-
[[e_i,[e_{i},e_{i-1}]_q]_{\q},e_{i+1}]_{q^2} =0$ from (9b).
If $i<m $  again from $Id(1)$ $z=0$. Hence
$[e_i,a_{i-1}^-]=0$, if $m\ne i$. 
\hfill\break
So far we have from (ii) and (iii) that
$$
[e_i,a_{i-1}^-]=0. \e(22)
$$
The rest of the proof is by induction. Assume $[e_i,a_j^-]=0$
for a certein $i>j,\;i\ne N$. Then from
(9a), (14a) and $Id(3)$ 
$[e_i,a_{j-1}^-]\sim [e_i,[e_{j-1},a_j^-]_{q_{j-1}}]=
[e_{j-1},[e_{i},a_j^-]]_{q_{j-1}}=0. $ Therefore
$[e_i,a_j^-]=0$ for any  $i>j,\;i\ne N$. Combining
the last with (i), one obtains  (18a).

\s\n
4. Eq.(18b) follows from the definition of $a_i^-$ and the
observation that
$(-1)^{(m-i-1)\l i+1 \r-(m-i)\l i \r}=(-1)^{\l i+1 \r}.$

\s\n
5. It remains to verify (18c). Since 
$a_i^-\sim [[e_i,e_{i+1}]_{q_i}, a_{i+2}^-]_{q_{i+1}}$,
$\[e_i,a_i^-\]_{\q_{i-1}}\sim
\[e_i,[[e_i,e_{i+1}]_{q_i},a_{i+2}^-]_{q_{i+1}}\]_{\q_{i-1}}$
(from (18a) and $Id(3)$)
$=[z,a_{i+2}^-]_{q_{i+1}}$, where 
$z=\[e_i,[e_i,e_{i+1}]_{q_i}\]_{\q_{i-1}}$ .
If $i>m$,   $z=[e_i,[e_i,e_{i+1}]_{q}]_{\q}=0$ (see (9b));
if $i<m$,   $z=[e_i,[e_i,e_{i+1}]_{\q}]_{q}=0$ again from  (9b);
if $i=m$,   $z=\{e_m,[e_m,e_{m+1}]_{q}\}_{q}=e_m^2e_{m+1}
-q^2e_{m+1}e_m^2 =0$, since according to (9a) $e_m^2=0$.
Hence $\[e_i,a_i^-\]_{\q_{i-1}}=0$.

\s\n
6. Eqs. (19) are proved in a similar way as Eqs. (18). This
completes the proof of the proposition. 

\smallskip\n
{\bf Proposition 2.} {\it The deformed Green operators (13) 
generate $\Uq$.}

\smallskip\n
{\it Proof.} The proof is an immediate consequence of the
relations:
$$
\eqalignno{
& \[a_{i}^-,a_{i+1}^+\]=2L_{i+1}e_{i}, \quad i=1,2,,\ldots,N-1,&(23a)\cr
& \[a_{i+1}^-,a_{i}^+\]=-2(-1)^{\l i+1 \r} f_{i} \L_{i+1},\quad
  \quad i=1,2,\ldots,N-1, &(23b)\cr
& \[a_i^-,a_i^+\]=-2{{L_i - \L_i}\over q-{\bar q}},\quad
  L_i=q^{H_i}, \quad\L_i=q^{-H_i}, \quad i=1,\ldots,N. & (23c) \cr
}
$$
These equations are proved by induction on $i$.

\s\n
1. The Serre relation (8c) together with the 
definitions of $a_N^\pm$, $L_N$ immediately yield
$
\[a_N^-,a_N^+\]=-2{{L_N - \L_N}\over q-{\bar q}}. 
$
From (18b) $a_{N-1}^- =\[e_{N-1},a_N^-\]_{q_{N-1}}
=[e_{N-1},a_N^-]_{q_{N-1}}$. 
Taking into account that $[e_{N-1},a_N^+]=0$ and $Id(3)$, one has
$
\[a_{N-1}^-,a_{N}^+\]=[a_{N-1}^-,a_{N}^+]
=[[e_{N-1},a_N^-]_{q_{N-1}},a_N^+]=
[e_{N-1}, [a_N^-,a_N^+]]_{q_{N-1}}=2L_Ne_{N-1},
$ i.e.,
$$
\[a_{N-1}^-,a_{N}^+\]=2L_Ne_{N-1}.\e(24a)
$$
In a similar way one shows that
$$
\[a_{N}^-,a_{N-1}^+\]=-2f_{N-1}\L_N.\e(24b)
$$
In order to compute $\[a_{N-1}^-,a_{N-1}^+\]$ set from (13b)
$a_{N-1}^+=-[a_N^+,f_{N-1}]_{\q_{N-1}}$. Since 
$n\ge 1, \; \q_{N-1}=\q$. Therefore
$\[a_{N-1}^-,a_{N-1}^+\]=-\[a_{N-1}^-,[a_N^+,f_{N-1}]_{\q}\]$.
Apply to the last supercommutator the identity $Id(2)$ with
$y=y=r=1$ and $x=s=t=\q$, namely
$$
\[A,[B,C]_{\q}\]=\[\[A,B\],C\]_{\q} +
(-1)^{deg(A)deg(B)}\[B,\[A,C\]\]_{\q},\e(25)
$$
where $A=a_{N-1}^-,\;\;B=a_N^+,\;\; C=f_{N-1}$. Then
\hfill\break
$\[a_{N-1}^-,a_{N-1}^+\]=-\[[a_{N-1}^-,a_N^+],f_{N-1}\]_{\q}-
[a_N^+,\[a_{N-1}^-,f_{N-1}\]]_{\q}= $(from (17) and (24a))
$=-\[2k_Ne_{N-1},f_{N-1}\]_{\q}-\[a_N^+,a_N^-\k_{N-1}\]_{\q}=
-2k_N\[e_{N-1},f_{N-1}\]+ [a_N^-,a_N^+]\k_{N-1}=
\hfill\break
-2k_N{{k_{N-1}-\k_{N-1}}\over{q-\q}}
-2{{k_{N}-\k_{N}}\over{q-\q}}\k_{N-1}=
-2{{k_{N}k_{N-1}-\k_{N}\k_{N-1}}\over{q-\q}}$, i.e.,
$$
\[a_{N-1}^-,a_{N-1}^+\]=-2{{L_{N-1}-\L_{N-1}}\over{q-\q}}.\e(26)
$$
From (24) and (26) we conclude that Eqs. (23) are fulfilled for
$i=N-1$.

\s\n
2. Assume Eqs. (23) hold for $i$ replaced by $i+1$:
$$
\eqalignno{
& \[a_{i+1}^-,a_{i+2}^+\]=2L_{i+2}e_{i+1}, & (27a) \cr
& \[a_{i+2}^-,a_{i+1}^+\]=-2(-1)^{\l i+2 \r} f_{i+1} \L_{i+2}, &(27b)\cr
& \[a_{i+1}^-,a_{i+1}^+\]=-2{{L_{i+1}-\L_{i+1}}\over q-{\bar q}}.&(27c) \cr
}
$$
We proceed to show that then Eqs. (23) hold too.  Set from (13b)
$a_i^-=(-1)^{\l i+1 \r}[e_i,a_{i+1}^-]_{q_i} $. Take into account that
according to (16) $[e_i,a_{i+1}^+]=0$ and $Id(3)$. Then
$\[a_i^-,a_{i+1}^+\]=(-1)^{\l i+1 \r} \[ [e_i,a_{i+1}^-]_{q_i},a_{i+1}^+\]
= (-1)^{\l i+1 \r}[ e_i,\[a_{i+1}^-,a_{i+1}^+\] ]_{q_i} $
(from (27a)$=2(\q-q)^{-1}(-1)^{\l i+1 \r}[e_i,L_{i+1}-\L_{i+1}]_{q_i}$
which after some rearrangement of the multiples finally yields
(23a). The verification of (23b) is similar. So far
we have derived from (27) that
$$
\eqalignno{
& \[a_{i}^-,a_{i+1}^+\]=2L_{i+1}e_{i}, \quad i=1,2,,\ldots,N-1,&(28a)\cr
& \[a_{i+1}^-,a_{i}^+\]=-2(-1)^{\l i+1 \r} f_{i} \L_{i+1},\quad
  \quad i=1,2,\ldots,N-1, &(28b)\cr
}
$$
Set in $\[a_i^-,a_i^+\]$ $\;\;a_i^+=-[a_{i+1}^+,f_i]_{\q_i}$. Either
$a_i^+$ or $f_i$ is an even element. Therefore applying 
again the identity (25),  one has
$
\[a_i^-,a_i^+\]=-\[a_i^-,[a_{i+1}^+,f_i]_{\q_i}\]=
-\[ \[a_i^-,a_{i+1}^+\],f_i\]_{\q_i}
-(-1)^{\l i \r \l i+1 \r}\[a_{i+1}^+,\[a_i^-,f_i\]\]_{\q_i}  
$
(from (28a) and (17)
$
=-2\[L_{i+1}e_i,f_i\]_{\q_i}
-(-1)^{\l i \r \l i+1 \r}\[a_{i+1}^+,a_{i+1}^-\k_i\]_{\q_i}=
-2L_{i+1}\[e_i,f_i\]+\[a_{i+1}^-,a_{i+1}^+\]\k_i=
-2{{L_i - \L_i}\over q-{\bar q}}$. Hence (23c) holds too.
From here and (28) we conclude that, if Eqs.(27) hold, then also
Eqs. (23) are fulfilled too. This completed the proof of the validity
of Eqs. (23).

From (13) and (23a,b) one obtains ($i=1,\ldots,N-1$):
$$
\eqalignno{
& h_i=H_i-H_{i+1}, \quad H_N=h_N, & (29a)\cr
& e_i={1\over 2}\L_{i+1}\[a_i^-,a_{i+1}^+\],
\quad e_N={1\over \sqrt{2}}a_N^-, & (29b) \cr
& f_i=-{1\over 2}(-)^{\l i+1 \r}\[a_{i+1}^-,a_{i}^+\]L_{i+1}
     ={1\over 2}\[a_i^+,a_{i+1}^-\]L_{i+1} ,
\quad f_N=-{1\over \sqrt{2}}a_N^+. & (29c) \cr
}
$$
Since the Chevalley elements generate $\Uq$, so do the Green
operators. This completes the proof.

\smallskip\n
{\bf Proposition 3.} {\it The Green generators 
$H_i,\; a_i^\pm,\; i=1,\ldots,N$ satisfy the following
relations 
\hfill\break
($i,j=1,\ldots,N,\;\;
\xi,\eta =\pm\; or\; \pm 1$ ) }:
$$
\eqalignno{
& [H_i,H_j]=0, & (30a) \cr
& [H_i,a_j^\pm]=\pm \delta_{ij}(-1)^{\l i \r} a_j^\pm, & (30b)\cr
& \[a_i^-,a_i^+\]=-2{L_i-{\bar L}_i\over q-{\bar q}},& (30c) \cr 
& [[a_{N-1}^\xi,a_N^\xi],a_{N}^\xi]_{\bar q}=0, & (30d) \cr
& \[\[a_i^{\eta},a_{i+ \xi}^{-\eta}\],  
a_j^{\eta}\]_{q^{-\xi (-1)^{\l i \r}\delta_{ij}}}
=2(\eta)^{\l j \r}\delta _{j,i + \xi}L_j^{-\xi \eta}a_i^{\eta}.&(30e)  \cr  
}
$$

\n
{\it Proof.}
The commutation relations (30a) are evident. (30b) follows from
the definitions of the Green generators, the Cartan relations
(8a,b) and the observation that
$$
\sum_{s=i}^N \sum_{r=j}^N \alpha_{sr}=-(-1)^{\l i \r}
\delta_{ij}. \e(31)
$$
The Eq. (30c) was derived in Proposition 2. The Eq.
$[[a_{N-1}^-,a_N^-],a_{N}^-]_{\bar q}=0$ is the same as
the Serre relation (9b), if one takes into account that
$[e_{N-1},e_N]_q\sim a_{N-1}^- $ and $e_N\sim a_N^-$.
Similarly one shows that (30d) with $\xi=+$ is the same as (10d).

The proof of (30c) is based on a case by case considerations
($\xi, \eta =\pm $).
To this end  one has to replace $e_i$ and $f_i$ in Eqs. (16)-(19)
with their expressions through the CAOs from (29). Using the relations
(which follow from (30b))
$$
\eqalignno{
& L_ia_j^\pm=a_j^\pm L_i,\quad \L_ia_j^\pm=a_j^\pm\L_i,
\quad i\ne j=1,\ldots,N, & (32a) \cr
& L_ia_i^\pm=q^{\pm (-1)^{\l i \r} } a_i^\pm L_i,\quad
\L_ia_i^\pm=q^{\mp (-1)^{\l i \r} } a_i^\pm \L_i,\quad i=1,\ldots,N, & (32b)
 \cr 
}
$$
after long, but simple calculations one verifies (30c).

\vskip 24pt

\n
{\bf 3. Description of osp$_{\bf q}$(2n+1/2m) via
Deformed Green Generators}

\bigskip
So far we have established that the Green generators
(13) satisfy the relations (30). Here we solve the inverse
problem: we show that the operators $H_i,\;
a_i^\pm,\;\;i=1,\ldots,N$ subject to the relations (30) provide
an alternative description of $\Uq$.

In Sect. 2 we have derived the relations (16)-(19)
from the definition (13) of the Green generators and the 
Cartan-Kac and the Serre relations, satisfied by the Chevalley
generators. Now as a first step we derive (16)-(19) on the ground
of Eqs. (30).

\smallskip\n
{\bf Proposition 4.} {\it The "mixed" relations (16)-(19)
follow from (29) and (30).}

\smallskip\n
{\it Proof.} Consider Eq. (16). Since $i\ne N$, from (29b)
$$
\[e_i,a_j^+\]={1\over 2}\L_{i+1}\[a_i^-,a_{i+1}^+\]a_j^+
-{1\over 2}a_j^+\L_{i+1}\[a_i^-,a_{i+1}^+\].\e(33)  
$$
\smallskip\n
(i) If $i+1<j$ or $i>j$, then $\L_{i+1}$ and $a_j^+$ commute (see (32)).
Therefore, using (30e), \hfill\break
$\[e_i,a_j^+\]={1\over 2}\L_{i+1}\[\[a_i^-,a_{i+1}^+\],a_j^+\]
=-{1\over 2}(-1)^{\l i+1 \r}\L_{i+1} \[\[a_{i+1}^+,a_{i}^-\],a_j^+\]=0 .$

\smallskip
\noindent
(ii) If $i=j$, again from (30e) $\[e_i,a_i^+\]
=-{1\over 2}(-1)^{\l i+1 \r}\L_{i+1} \[\[a_{i+1}^+,a_{i}^-\],a_i^+\]
=-(-1)^{\l i+1 \r}\L_{i+1}L_ia_{i+1}^+ \hfill\break
=-(-1)^{\l i+1 \r}k_ia_{i+1}^+$.

\smallskip\n
(iii) If $i+1=j$, from (33) and taking into account (32b)
\hfill\break
$\[e_i,a_{i+1}^+\]={1\over 2}\L_{i+1}(\[a_i^-,a_{i+1}^+\]a_{i+1}^+
-q^{(-1)^{\l i+1 \r} }a_{i+1}^+ \[a_i^-,a_{i+1}^+\])=
{1\over 2}\L_{i+1}\[\[a_i^-,a_{i+1}^+\],a_{i+1}^+\]_{q^{(-1)^{\l
i+1 \r} }} \hfill\break
={1\over 2}\L_{i+1}(-1)^{\l i \r \l i+1 \r}
\[\[a_{i+1}^+,a_i^-\],a_{i+1}^+\]_{q^{(-1)^{\l i+1 \r} }}=0.
$
The unification of (i)-(iii) yields (16). 

\n
The remaining equalities (17)-(19) are proved in a similar way.

\smallskip\n
{\bf Proposition 5.} {\it The Cartan-Kac relations are a
consequence of the relations (30).}

\smallskip\n
{\it Proof.} The first two equations (8a) and (8b) are easily
verified. We proceed to prove (8c).

\smallskip\n
1. The case $i=j$. If $i=N$, then (8c) is the same as (30c).
Let $i<N$. From (29c) and the graded Leibnitz rule 
\hfill\break
$\[e_i,f_i\]={1\over 2}\[e_i,\[a_i^+,a_{i+1}^-L_{i+1}\]\]=
{1\over 2}\[\[e_i,a_i^+\],a_{i+1}^-L_{i+1}\]
+(-1)^{(\l i \r + \l i+1 \r)\l i \r}\[a_i^+,\[e_i,a_{i+1}^-L_{i+1}\]\].$
\hfill\break
Insert above
\hfill\break
$\[e_i,a_i^+\]=-(-1)^{\l i+1 \r}k_i a_{i+1}^+ $ 
and 
$\[e_i,a_{i+1}^-L_{i+1}\] =\[e_i,a_{i+1}^-\]_{q_i}L_{i+1}$
(from (18b))=$(-1)^{\l i+1 \r}a_i^-L_{i+1}$.
\hfill\break
After some rearrangement of the multiples one obtains:
\hfill\break
$\[e_i,f_i\]={1\over 2}\[a_{i+1}^-,a_{i+1}^+\]L_i
-{1\over 2}\[a_{i}^-,a_{i}^+\]L_{i+1}$ (from (23c)
$={{k_i-{\bar k}_i}\over{q-{\bar q}}} $.
\hfill\break
Hence
$$
\[e_i,f_i\]={{k_i-{\bar k}_i}\over{q-{\bar q}}}
\quad i=1,\ldots,N.\e(34)
$$

\smallskip\n
2. The case $i\ne j$. Eq.
(8c) is easily verified for $i=N$ or $j=N$. We consider $i\ne j
\ne N.$ From (29c)
\hfill\break
$\[e_i,f_j\]={1\over 2}\[e_i,\[a_j^+,a_{j+1}^-L_{j+1}\]\]=
{1\over 2}\[\[e_i,a_j^+\],a_{j+1}^-L_{j+1}\]
+(-1)^{(\l i \r + \l i+1 \r)\l j \r}\[a_j^+,\[e_i,a_{j+1}^-L_{j+1}\]\].$
\hfill\break
The first term in the r.h.s. cancels out, since 
$\[e_i,a_j^+\]=0$ according to (16). From the second term evaluate
only the internal supercommutator $A=\[e_i,a_{j+1}^-L_{j+1}\]$.

\smallskip\n
(i) If $i<j$ or $i>j+1$, $A=\[e_i,a_{j+1}^-\]L_{j+1}=0$
according to (18a);

\smallskip\n
(ii) If $i=j+1$, then $L_ie_i=\q_{i-1}e_iL_i$.
Therefore $A=\[e_i,a_{i}^-\]_{\q_{i-1}} L_{i}=0$ 
according to (18c).

\smallskip\n
Hence 
$\[e_i,f_j\]=0$, if $i\ne j=1,\ldots,N$. The latter
together with (34) shows that also the last Cartan-Kac relation
(8c) is fulfilled. This completes the proof.

\smallskip\n
{\bf Proposition 6.} {\it The Serre relations (9) and (10) are a
consequence of the relations (30).}

\smallskip\n
{\it Proof.} 

\n 1. First we prove that $[e_i,e_j]=0$, if $|i-j|>1$. Assume for
definiteness that $i+1<j.$

\smallskip
\n (i) Let $i+1<j=N$. From (29b) and the observation that $\L_{i+1}$
commutes with $a_N^-$ one has
\hfill\break
$[e_i,e_N]\sim [\L_{i+1}\[a_i^-,a_{i+1}^+\],a_N^-]=
 \L_{i+1}[\[a_i^-,a_{i+1}^+\],a_N^-]$ (from (30e))$=0$.

\smallskip
\n (ii) Let $i+1<j<N$. From (29b)
$[e_i,e_j]\sim [\L_{i+1}\[a_i^-,a_{i+1}^+\],
 \L_{j+1}\[a_j^-,a_{j+1}^+\]]$ ($\L_{j+1}$ commutes with
$ a_i^-$ and $a_{i+1}^+$, $\L_{i+1}$ commutes with
$ a_j^-$ and $a_{j+1}^+$)
$= \L_{i+1} \L_{j+1} \[\[a_i^-,a_{i+1}^+\],\[a_j^-,a_{j+1}^+\]\]
\hfill\break
=\L_{i+1} \L_{j+1} \[\[\[a_i^-,a_{i+1}^+\],a_j^-\],a_{j+1}^+\]\]
+(-1)^{(\l i \r + \l i+1 \r)\l j \r}\L_{i+1} \L_{j+1}
\[a_j^-,\[\[a_i^-,a_{i+1}^+\],a_{j+1}^+\]\]=0$, since
\hfill\break
$\[\[a_i^-,a_{i+1}^+\],a_j^-\]=0$ and
$\[\[a_i^-,a_{i+1}^+\],a_{j+1}^+\]=0$ according to (30e).

\s\n
2. The proof of  $[f_i,f_j]=0$, if $|i-j|>1$ is similar.

\smallskip\n
3. Proof of $[e_i,[e_i,e_{i+1}]_{q'}]_{\q'}=0, \quad
  i\neq m,\;\; i\neq N$ and $q'=q$ or $q'=\q$.

\smallskip\n
We choose $q'=q_{i-1}=q^{(-1)^{\l i \r}}$. Therefore
the relation to be proved is
$$
[e_i,[e_i,e_{i+1}]_{q_{i-1}}]_{\q_{i-1}}=0, \quad
  i\neq m, \;\; i\neq N. \e(35)
$$
As a preliminary step compute $[e_i,e_{i+1}]_{q_{i-1}}$ (see (29b))
$=  {1\over 4}[\L_{i+1}\[a_i^-,a_{i+1}^+\],
     \L_{i+2}\[a_{i+1}^-,a_{i+2}^+\]]_{q_{i-1}}$. 
\hfill\break	 
From (32) one has
\hfill\break
 $\[a_i^-,a_{i+1}^+\]\L_{i+2}=\L_{i+2}\[a_i^-,a_{i+1}^+\]$ and
 $\[a_{i+1}^-,a_{i+2}^+\]\L_{i+1}=
 q^{-(-1)^{\l i+1 \r}}\L_{i+1}\[a_{i+1}^-,a_{i+2}^+\]
= \q_i\L_{i+1}\[a_{i+1}^-,a_{i+2}^+\]$.
\hfill\break
Therefore, 
$[e_i,e_{i+1}]_{q_{i-1}}={1\over 4}\L_{i+1}\L_{i+2}
\left(\[a_i^-,a_{i+1}^+\]\[a_{i+1}^-,a_{i+2}^+\]-
q_{i-1}\q_i\[a_{i+1}^-,a_{i+2}^+\]\[a_i^-,a_{i+1}^+\]\right).$
\hfill\break
Since $i\ne m$, $q_{i-1}\q_i=1$. Thus,
$[e_i,e_{i+1}]_{q_{i-1}}={1\over 4}\L_{i+1}\L_{i+2}
[\[a_i^-,a_{i+1}^+\],\[a_{i+1}^-,a_{i+2}^+\]]$
\hfill\break
$={1\over 4}\L_{i+1}\L_{i+2}
\[\[\[a_i^-,a_{i+1}^+\],a_{i+1}^-\],a_{i+2}^+\]
+{1\over 4}\L_{i+1}\L_{i+2}(-1)^{(\l i \r + \l i+1 \r)\l i+1 \r}
\[a_{i+1}^-,\[\[a_i^-,a_{i+1}^+\],a_{i+2}^+\]\].$
\hfill\break
The second term in the r.h.s. is zero, since from (30c)
$\[\[a_{i+1}^+,a_i^-\],a_{i+2}^+\]=0$. Again from
(30c) 
\hfill\break
$\[\[a_i^-,a_{i+1}^+\],a_{i+1}^-\]=
2(-1)^{\l i+1 \r} L_{i+1}a_i^-.$ 
Therefore
$$
[e_i,e_{i+1}]_{q_{i-1}}={1\over 2}(-1)^{\l i+1 \r}\L_{i+2}
\[a_i^-,a_{i+2}^+\].\e(36)
$$
Insert $e_i$ from (29b) and $[e_i,e_{i+1}]_{q_{i-1}}$ from (36) 
in the l.h.s. of (35):
\hfill\break
$[e_i,[e_i,e_{i+1}]_{q_{i-1}}]_{\q_{i-1}}=
[{1\over 2}\L_{i+1}\[a_i^-,a_{i+1}^+\],
{1\over 2}(-1)^{\l i+1 \r}\L_{i+2}\[a_i^-,a_{i+2}^+\]]_{\q_{i-1}}
\hfill\break
={1\over 4}(-1)^{\l i+1 \r}\L_{i+1}\L_{i+2} [ \[a_i^-,a_{i+1}^+\],
\[a_i^-,a_{i+2}^+\]]_{\q_{i-1}}$ (use the circumstance that
$[\[a_i^-,a_{i+1}^+\],a_{i+2}^+]=0$ and $Id(3)$)
$={1\over 4}(-1)^{\l i+1 \r}\L_{i+1}\L_{i+2} 
\[ [ \[a_i^-,a_{i+1}^+\],a_i^-]_{\q_{i-1}},a_{i+2}^+\]=0$,
since 
$[ \[a_i^-,a_{i+1}^+\],a_i^-]_{\q_{i-1}}=
[ \[a_i^-,a_{i+1}^+\],a_i^-]_{q^{-(-1)^{\l i \r}}}=0$ according
to (30e). Hence the Serre relation (35) holds.

\smallskip\n
4. The proof of $[e_i,[e_i,e_{i-1}]_{q'}]_{\q'}=0, \quad
   i\neq 1,\;\;  i\neq m,\;\; i\neq N$ and $q'=q$ or $q'=\q$
   is similar. For $q'$ one has to take 
   $q'=\q_{i-1}=q^{-(-1)^{\l i \r}}$.
  
\s\n
5. The proof of the Serre relations (10b) is similar as for
   (9b).

\smallskip\n
6. Proof of $e_m^2=0$ (i.e., of (9a) for $i=j=m$).

\smallskip\n
$e_m^2\sim \[e_m,e_m\]_{q^2}$ (use (29b))
$\sim \[\L_{m+1}\[a_m^-,a_{m+1}^+\],\L_{m+1}\[a_m^-,a_{m+1}^+\]\]_{q^2}
=q\L_{m+1}^2\[\[a_m^-,a_{m+1}^+\],\[a_m^-,a_{m+1}^+\]\]_{q^2}$.
\hfill\break
In order to evaluate
$\[\[a_m^-,a_{m+1}^+\],\[a_m^-,a_{m+1}^+\]\]_{q^2}$
set $A=\[a_m^-,a_{m+1}^+\]$, $B=a_{m}^-$, $C=a_{m+1}^+$.
Note that $deg(A)=deg(B)={\bar 1}$, $deg(C)={\bar 0}$
and use the identity $Id(2)$ with 
$x=1,\; y=q^2,\; z=r=t=q, \; s=\q$. It yields
\hfill\break
$\[\[a_m^-,a_{m+1}^+\],\[a_m^-,a_{m+1}^+\]\]_{q^2}=
\[\[\[a_m^-,a_{m+1}^+\],a_m^-\]_q,a_{m+1}^+\]_{q}-
q\[a_m^-,\[\[a_m^-,a_{m+1}^+\],a_{m+1}^+\]_{q}\]_{\q}$=0,
\hfill\break
since, as it follows from (30e),
$\[\[a_m^-,a_{m+1}^+\],a_m^-\]_q=0$ and
$\[\[a_m^-,a_{m+1}^+\],a_{m+1}^+\]_{q}=0 $.

\s\n 7. The proof of $f_m^2=0$ is similar.

\s\n
8. Proof of $\{[e_{m},e_{m-1}]_{q},[e_m,e_{m+1}]_{\bar q}\}=0$.

\s\n 
So far we have proved the validity of the Cartan-Kac relations (8)
and of the Serre relations (9a,b) and (10a,b). Therefore we
can refer to them. In particular Eqs. (15) hold.
Using (15a), write 
$a_{m-1}^-=-[[e_{m-1},[e_{m},e_{m+1}]_{q}]_{\q},a_{m+2}^-]_{q}$.
According to (18a), which was proved in Proposition 4,
 $\[e_m,a_{m-1}^-\]=0$. Therefore,
$0=\[e_m,a_{m-1}^-\]
=-\[e_m,[[e_{m-1},[e_{m},e_{m+1}]_{q}]_{\q},a_{m+2}^-]_{q}\]$
and since, again from (18a), $[e_m,a_{m+2}^-]=0$, applying
$Id(3)$, we have
$\[e_m,a_{m-1}^-\]=-[y,a_{m+2}^-]_q$, where
$y=\[e_{m},[e_{m-1},[e_{m},e_{m+1}]_{q}]_{\q}\]$, which can be written
also as
$$
y=\[e_{m},[[e_{m-1},e_{m}]_{\q},e_{m+1}]_{q}\]\e(37)
$$  
and
$$
[y,a_{m+2}^-]_q=0. \e(38)
$$   
From (13b), (37) and (8c) one immediately concludes that
$[y,a_{m+2}^+]=0.$ Therefore, applying $Id(3)$, one has
$0=[[y,a_{m+2}^-]_q,a_{m+2}^+]=[y,[a_{m+2}^-,a_{m+2}^+]]_q $
(use (30c))$=-2(q-\q)^{-1}[y,L_{m+2}-\L_{m+2}]_q$, which
after, pushing $L_{m+2}$ and $\L_{m+2}$ to the right, yields
$(1-q^2)yL_{m+2}=0$. Hence
$$
y=\[e_{m},[[e_{m-1},e_{m}]_{\q},e_{m+1}]_{q}\]=0. \e(39)
$$
Set in (39) $A=e_m,\;\; C=[e_{m-1},e_{m}]_{\q},\;\; B=e_{m+1}$
and use the following identity, which follows from $Id(2)$:
If $B$ is an even element, then
$$
\[A,[C,B]_q\]=-q\[ [A,B]_{\q},C\]-[B,\[A,C\]_q ].\e(40)
$$
This yields $y=-q\[ [e_{m},e_{m+1}]_{\q},[e_{m-1},e_{m}]_{\q}\]-
[e_{m+1},\[e_{m}, [e_{m-1},e_{m}]_{\q}\]_{q}]=
\{ [e_{m},e_{m-1}]_{q},[e_{m},e_{m+1}]_{\q}\}
-[e_{m+1},\{e_{m}, [e_{m-1},e_{m}]_{\q}\}_{q}]=0.$
The second term in the r.h.s. above is zero, since
\hfill\break
$ \{e_{m}, [e_{m-1},e_{m}]_{\q}\}_{q}
=qe_{m-1}e_m^2-\q e_m^2e_{m-1}$ and $e_m^2=0$.
Therefore,
$$
\{ [e_{m},e_{m-1}]_{q},[e_{m},e_{m+1}]_{\q}\}=0,
$$
which proves (9c).

\s\n
9. The proof of (10c) is similar.

\s\n
10. Proof of $[e_N,[e_N,[e_N,e_{N-1}]_{\bar q}]]_q=0$.

\s\n From (29b) $e_{N-1}={1\over 2}\L_{N}[a_{N-1}^-,a_{N}^+]$ and
$e_N={1\over \sqrt{2}}a_N^-$. Therefore,
$[e_N,e_{N-1}]_{\q}=
[{1\over \sqrt{2}}a_N^-,{1\over 2}\L_{N}[a_{N-1}^-,a_{N}^+]]_{\q}=
{1\over{2\sqrt 2}}\q\L_n[a_N^-,[a_{N-1}^-,a_N^+]]$, which,
applying (30e), yields:
$$
[e_N,e_{N-1}]_{\q}=-{\q\over \sqrt 2}a_{N-1}^-. \e(41)
$$
Therefore,
$[e_N,[e_N,[e_N,e_{N-1}]_{\bar q}]]_q
=-{\q\over 2\sqrt 2}[a_N^-,[a_N^-,a_{N-1}^-]]_{q}
=-{1\over 2\sqrt 2}[[a_{N-1}^-,a_N^-],a_{N}^-]]_{\q}=0$,
according to (30d). Hence, (9d) holds.

\s\n
11. The proof of (10d) is similar.

\s\n
This completes the proof of Proposition 6.

The relations (29), written in the form
($i=1,\ldots,N-1$)
$$
\eqalignno{
& h_i=H_i-H_{i+1}, \quad H_N=h_N, & (42a)\cr
& e_i={1\over 2}q^{-H_{i+1}}\[a_i^-,a_{i+1}^+\],
\quad e_N={1\over \sqrt{2}}a_N^-, & (42b) \cr
& f_i={1\over 2}\[a_i^+,a_{i+1}^-\]q^{H_{i+1}} ,
\quad f_N=-{1\over \sqrt{2}}a_N^+, & (42c) \cr
}
$$
indicate  that the Chevalley elements are functions of the Green
generators. More precisely, $h_i,e_i,f_i$ are in the closure of
the subalgebra of all polynomials of the Green operators over
$\Ch$.

\bigskip
So far we were considering $\Uq$ as a topologically free module
over the ring $\Ch$ of the formal power series over an
indeterminate $h$.  Due to this, for instance, $q^{H_i}$ is a
well defined element from $\U_q$. It is important to note however
that all our considerations remain true, if one goes to the
factor algebra $U_{h_c}[osp(2n+1/2m)]$, replacing $h$ by a
complex number $h_c$, such that $h_c\notin i\pi \Q$ ($\Q$ - all
rational numbers), namely considering $q$ to be a number, which
is not a root of 1. Then in the limit $h_c \rightarrow 0$ the
deformed Green generators become ordinary parabosons and
parafermoins.  This is the justification to call the operators
(13) deformed Green generators, and the corresponding to them
statistics - quantum deformation of the parastatistics. 
We conclude the paper, formulating our main result as a
theorem.

\bigskip\n
{\bf Theorem.} {\it $\Uq$ is a topologically free $\Ch$ module
and a unital algebra with generators $H_i,\;
a_i^\pm,\;\;i=1,\ldots,N$ and relations (30). The generators
consist of $m$ pairs of deformed parabosons and $n$ pairs of deformed
parafermions.  }

\bigskip
This theorem established a link between the quantum
groups in the sense of Drinfeld-Jimbo$^{3,4}$ and the quantum
statistics in the sense of Green.$^2$

\vskip 18pt
\noindent
{\it Acknowledgments.}
I am grateful to Prof. Randjbar-Daemi for the kind hospitality at
the High Energy Section of ICTP.  This work was supported by the
Grant $\Phi-416$ of the Bulgarian Foundation for Scientific
Research.

\bigskip\n
{\bf References}

\vskip 12pt
{\settabs \+  $^{11}\;\;$ & I. Patera, T. D. Palev, Theoretical
   interpretation of the experiments on the elastic \cr

\+ $^1$ & T.D. Palev, J. Phys. A~: Math. Gen. {\bf 29}, L171 (1996).\cr

\+ $^2$  & H.S. Green, Phys. Rev. {\bf 90}, 270 (1953).  \cr

\+ $^3$  & V. Drinfeld, {\it Quantum groups.} ICM proceedings,
            Berkeley 798 (1986).\cr
\+ $^4$  & M. Jimbo, Lett. Math. Phys. {\bf 11}, 247 (1986).\cr

\+ $^5$  & S. Kamefuchi and Y. Takahashi, 
           Nucl. Phys. {\bf 36}, 177 (1962). \cr

\+ $^6$  & C. Ryan  and E.C.G. Sudarshan, 
           Nucl. Phys. {\bf 47}, 207 (1963). \cr

\+ $^7$  & M. Omote, Y. Ohnuki and S. Kamefuchi,
           Progr. Theor. Phys. {\bf 56} 1948 (1976).\cr

\+ $^8$  & A.Ch. Ganchev and T.D. Palev, 
           J. Math. Phys. {\bf 21}, 797 (1980).\cr

\+ $^9$  & O.W. Greenberg and A.M. Messiah,
           Phys. Rev. {\bf 138}, B1155 (1965).\cr

\+ $^{10}$ & T.D. Palev,
            J. Math. Phys. {\bf 23}, 1100 (1982).\cr

\+ $^{11}$ & A.J. Bracken  and H.S. Green,  Journ. Math. Phys. 
             {\bf 14} 1784 (1973) \cr
\+ $^{12}$ & T.D Palev,  Ann. Inst. Henri Poincare 
             {\bf XXIII} 49 (1975).\cr

\+ $^{13}$ & V.G. Kac, 
            Lecture Notes in Math. {\bf 676}, 597 (Springer, 1979).\cr

\+ $^{14}$ & T.D. Palev, Lie algebraical aspects of the quantum
               statistics. Thesis, Institute of Nuclear Research \cr 
\+ &           and Nuclear Energy (1976),  Sofia. \cr

\+ $^{15}$ & T.D. Palev, {\it Lie algebraical aspects of quantum statistics. 
            Unitary quantization (A-quantization)}.\cr
\+        & Preprint JINR E17-10550 (1977) and hep-th/9705032.\cr

\+ $^{16}$ & T.D. Palev, Rep. Math. Phys. {\bf 31}, 241 (1992).\cr

\+ $^{17}$ & T.D. Palev, Czech. Journ. Phys. {\bf B32}, 680 (1982).\cr

\+ $^{18}$ & T.D. Palev, J. Math. Phys. {\bf 23}, 1778 (1982).\cr

\+ $^{19}$ & T.D. Palev and N.I. Stoilova, 
            J. Math. Phys. {\bf 38}, 2506 (1997) and 
            hep-th/9606011. \cr

\+ $^{20}$ & F.D.M. Haldane, Phys. Rev. Lett. {\bf 67}, 937 (1991).\cr

\+ $^{21}$ & T.D. Palev and N.I. Stoilova, J. Phys. A~: Math. Gen. 
            {\bf 27}, 7387 (1994) and hep-th/9405125.\cr

\+ $^{22}$ & A.J. Macfarlane, J. Phys. A~: Math. Gen.
    	    {\bf 22}, 4581 (1989).  \cr

\+ $^{23}$ & L.C. Biedenharn,  J. Phys. A~: Math. Gen. {\bf 22}, L873. \cr

\+ $^{24}$ & C.P. Sun and H.C.Fu, J. Phys. A~: Math. Gen. 
	        {\bf 22},  L983.  \cr

\+ $^{25}$ & T.D. Palev, Lett. Math. Phys. {\bf 31}, 151 (1994)
            and hep-th/9311163.\cr

\+ $^{26}$ & E. Celeghini, T.D. Palev,  and M. Tarlini,
            Mod. Phys. Lett. {\bf B5}, 187 (1991).\cr

\+ $^{27}$ & T.D. Palev, J. Phys. A~: Math. Gen. {\bf 26}, L1111 (1993)
            and hep-th/9306016. \cr

\+ $^{28}$ & L.K. Hadjiivanov, J. Math. Phys. {\bf 34}, 5476 (1993). \cr
		 
\+ $^{29}$ & T.D. Palev and J. Van der Jeugt,
             J. Phys. A~: Math. Gen. {\bf 28}, 2605 (1995)
             and $q$-alg/9501020 .\cr
   
\+ $^{30}$ & M. Chaichian  and P. Kulish,  Phys. Lett. {\bf 234B} 72 (1990).\cr

\+ $^{31}$ & A.J. Bracken, M.D. Gould and R.B. Zhang, Mod. Phys. Lett. A
         {\bf 5} 831 (1990). \cr

\+ $^{32}$ & R. Floreanini, V.P. Spiridonov  and L. Vinet,
             Comm. Math. Phys. {\bf 137} 149 (1990).\cr

\+ $^{33}$ & S.M. Khoroshkin  and V.N. Tolstoy,
             Comm. Math. Phys. {\bf 141} 599 (1991). \cr

\+ $^{34}$ & R. Floreanini, D.A. Leites and L. Vinet, 		
             Lett. Math. Phys. {\bf 23} 127 (1991).\cr

\+ $^{35}$ & M. Scheunert, Lett. Math. Phys. {\bf 24 }, 173 (1992). \cr

\end